\documentclass[final]{siamltex}

\usepackage{amsmath,color}
\usepackage{epsfig}
\usepackage{verbatim}
\usepackage{amssymb}
\usepackage{multirow}
\usepackage{subfigure}
\usepackage[linesnumbered,algoruled,lined]{algorithm2e}

\newcommand{\lra}{{\longrightarrow}}
\newcommand{\eproof}{\hfill\rule{2.2mm}{3.0mm}}

\newcommand{\Proof}{\noindent {\bf Proof.~~}}

\newcommand{\R}{{\mathbb R}}

\newcommand{\C}{{\mathbb C}}

\newcommand{\ep}{\varepsilon}

\renewcommand{\eqref}[1]{(\ref{#1})}
\newcommand{\inner}[1]{\langle #1 \rangle}

\newcommand{\mhsp}{\hspace{2em}}
\newcommand{\FT}[1]{\widehat{#1}}
\newcommand{\ul}[1]{\underline{#1}}

\newcommand{\A}{{\mathcal A}}

\newcommand{\tr}{{\rm tr}}
\newcommand{\argmin}{{\rm argmin}}

\newcommand{\vx}{{\mathbf x}}
\newcommand{\vy}{{\mathbf y}}
\newcommand{\vz}{{\mathbf z}}

\newcommand{\va}{{\mathbf a}}
\newcommand{\vb}{{\mathbf b}}

\newcommand{\vv}{{\mathbf v}}

\newcommand{\vf}{{\mathbf f}}
\newcommand{\vg}{{\mathbf g}}

\renewcommand{\H}{{\mathbb H}}
\newcommand{\HH}{{\mathbf H}}
\newcommand{\qHH}{\overline{\HH}}
\newcommand{\MM}{\mathbf M}

\newcommand{\F}{{\mathcal F}}

\newtheorem{prop}{Proposition}[section]
\newtheorem{lem}[prop]{Lemma}

\newtheorem{theo}[prop]{Theorem}

\renewcommand{\theequation}{\thesection.\arabic{equation}}
\renewcommand{\thefigure}{\arabic{figure}}

\title{Fast Rank One Alternating Minimization Algorithm for Phase Retrieval}

\author{Jian-Feng Cai\thanks{Department of Mathematics, The Hong Kong University of Science and Technology,
Clear Water Bay, Kowloon, Hong Kong. JFC was supported in part by HKRGC grant 16300616. ({\tt jfcai@ust.hk}).}
        \and Haixia Liu\thanks{Department of Mathematics, The Hong Kong University of Science and Technology,
Clear Water Bay, Kowloon, Hong Kong. ({\tt mahxliu@ust.hk})}
\and Yang Wang
\thanks{Department of Mathematics, The Hong Kong University of Science and Technology,
Clear Water Bay, Kowloon, Hong Kong. Part of the research work of YW was completed while the author was at Department of Mathematics, Michigan State University. This research was supported in part by the National Science Foundation grant DMS-1043032 and AFOSR grant FA9550-12-1-0455 and HKRGC grant 16306415. ({\tt yangwang@ust.hk}).}}

\begin{document}

\maketitle

\begin{abstract}
The {\it phase retrieval} problem is a fundamental problem in many fields, which is appealing for investigation. It is to recover the signal vector $\tilde{\vx}\in\mathbb{C}^d$ from a set of $N$ measurements $b_n=|\vf^*_n\tilde{\vx}|^2,\ n=1,\cdots, N$, where $\{\vf_n\}_{n=1}^N$ forms a frame of $\mathbb{C}^d$. 
Existing algorithms usually use a least squares fitting to the measurements, yielding a quartic polynomial minimization. 
In this paper, we employ a new strategy by splitting the variables, and we solve a bi-variate optimization problem that is quadratic in each of the variables.
An alternating gradient descent algorithm is proposed, and its convergence for any initialization is provided.
Since a larger step size is allowed due to the smaller Hessian, the alternating gradient descent algorithm converges faster than the gradient descent algorithm (known as the Wirtinger flow algorithm) applied to the quartic objective without splitting the variables. Numerical results illustrate that our proposed algorithm needs less iterations than Wirtinger flow to achieve the same accuracy.
\end{abstract}


\pagestyle{myheadings}
\thispagestyle{plain}
\markboth{Fast Rank One Alternating Minimization Algorithm for Phase Retrieval}{Jian-Feng Cai, Haixia Liu, Yang Wang}
\section{Introduction}

Let $f(x) \in L^2(\R^d)$. It is well known that the map $f \mapsto \FT f$, where $\FT f$ denotes the Fourier transform of $f$, is an isometry in $L^2(\R^d)$ and hence $f$ can be uniquely reconstructed from $\FT f$. In many applications such as X-ray crystallography, however, we can only measure the magnitude $|\FT f|$ of the Fourier transform. This raises the following question: Is it still possible to reconstruct $f$ from $|\FT f|$? This is the classic {\em phase retrieval} problem.

The phase retrieval problem has a natural generalization to finite dimensional Hilbert spaces. Such an extension has important applications in imaging, optics, communication, audio signal processing and more \cite{chai2010array,harrison1993phase,heinosaari2013quantum,millane1990phase,walther1963question}. It is in this finite Hilbert space setting that phase retrieval has become one of the growing areas of research in recent years. 

Let $\HH$ be a (real or complex) Hilbert space of finite dimension. Without loss of generality we identify $\HH$ with $\H^d$ where $\H=\R$ or $\H=\C$. A set of elements $\F=\{\vf_n\}$ in $\HH$ is called a {\em frame} if it spans $\HH$. Given this frame any vector $\vx\in {\HH}$ can be reconstructed from the inner products $\{\inner{\vx,\vf_n}\}$. Often it is convenient to identify the frame $\F$ with the corresponding {\em frame matrix} $F=[\vf_1,\vf_2, \dots, \vf_N]$. The phase retrieval problem in $\HH$ is: 

\medskip
\noindent
{\bf The Phase Retrieval Problem.}~~{\em Let $\F=\{\vf_n\}$ be a frame in $\HH$. Can we reconstruct any $\vx\in\HH$ up to a unimodular scalar from $\{|\inner{\vx,\vf_n}|\}$, and if so, how?}

\medskip

$\F$ is said to be {\em phase retrievable (PR)} if the answer is affirmative. There is an alternative formulation. Consider the equivalence relation $\sim$ on $\HH$: $\vx_1 \sim \vx_2$ if there is a constant
$b\in \H$ with $|b|=1$ such that $\vx_1=b\vx_2$. Let $\qHH :=\HH/\sim$. We shall use
$\ul\vx$ to denote the equivalent class containing $\vx$. For any given frame
$\F=\{\vf_n: 1\leq n \leq N\}$ in $\HH$ define the map $\MM_\F: \qHH \lra \R_+^N$ by
\begin{equation}  \label{1.1}
      \MM_\F(\ul\vx) = [|\inner{\vx,\vf_1}|^2, \dots, |\inner{\vx,\vf_N}|^2]^T.
\end{equation}
The phase retrieval problem asks whether a $\ul\vx\in\qHH$ is uniquely determined
by $\MM_\F(\ul\vx)$, i.e. whether $\MM_\F$ is {\em injective} on $\qHH$.

Many challenging and fundamental problems in phase retrieval remain open. For example, for phase retrieval in $\C^d$ it is still unknown what is the minimal number of vectors needed for a set of vectors $\F$ to be phase retrievable. A challenging problem of very practical importance is the computational efficiency of phase retrieval algorithms. So far the existing phase retrieval algorithms can be loosely divided into four categories: (A) Using frames with very large $N$, in the order of $N \geq O(d^2)$, (B) Convex relaxation algorithms using random frames, (C) Non-convex optimization with a quartic objective with random frames,  and (D) Constructing special frames $\F$ that allow for fast and robust phase retrieval reconstruction of $\vx$.

The first category is based on the fact that each $|\inner{\vx,\vf_n}|^2$ is a linear combination of monomials $x_i^*x_j$. The reconstruction of $\vx$ can be attained by solving for these monomials, provided that there are enough equations, i.e. $N$ is large enough. We will need $N \geq \frac{1}{2}d(d+1)$ in the real case and $N \geq d^2$ in the complex case. The reconstruction then becomes solving a system of linear equations if we treat all monomials as independent variables. One can also obtain robustness results under such framework. The weakness of this approach is that when $d$ is large the number of variables and the number of measurements needed will explode, making it generally impractical and slow. Several constructions for special frames were designed (e.g. \cite{balan2009painless}) with which one can compute $\vx$ efficiently (``painless reconstruction''). But this doesn't reduce the number of required measurements.

The second category of methods employ convex relaxation techniques like those of compressive sensing. By considering $X=\vx\vx^*$ we can rewrite the map $\MM_\F$ as
\begin{equation}  \label{1.2}
      \MM_\F(X) = [\vx^*A_1\vx, \dots, \vx^*A_N\vx]^T
                = [\tr(A_1X), \dots, \tr(A_NX)]^T
\end{equation}
where $A_n = \vf_n\vf^*_n$, which is a linear map from $\C^{d\times d}$ to $\R^N$. The original problem is now a linear equation $\MM_\F(X) = \vb$ subject to the constraints $X \geq 0$ and has rank 1. This type of problems is not convex and cannot be solved efficiently in general. However, it was shown in \cite{candes2013phaselift} that with high probability for random frames with $N \geq O(d\log d)$ the original problem is equivalent to the convex problem of solving for 
\begin{equation} \label{1.3}
   \argmin_X \tr(X) \mhsp\mbox{subject to} \mhsp X\geq 0, ~~\MM_\F(X) = \vb.
\end{equation}
The bound was later improved to $N \geq O(d)$ \cite{candes2014solving}. The robustness of the method was also proved. This convex relaxation method, called {\em PhaseLift},  solves (\ref{2.3})  using semi-definite programming. It can be done with reasonable efficiency for small $d$ (typically up to about $d=1000$ on a PC). Several refinements and variations of PhaseLift have also being proposed, e.g. PhaseCut and MaxCut \cite{candes2015phase1,waldspurger2015phase} for Fourier  measurements with random masks. But all employ semi-definite programming, which for larger $d$ becomes slow and impractical. 

The third category of methods consider the non-linear least square fittings to the measurements $|\langle\vx,\vf_n\rangle|^2$, $n=1,\ldots,N$, and solve
\begin{equation}\label{1.4}
\min_{\vx\in\mathbb{C}^d}\sum_{n=1}^{N}\left(|\langle\vx,\vf_n\rangle|^2-b_n\right)^2.
\end{equation}
Since the objective is a quartic polynomial, \eqref{1.4} is a smooth non-convex optimization. Such an optimization is very difficult to solve in general because the non-convex objective may have numerous local minima. However, when random frames are used, \eqref{1.4} is not as difficult as it appears. In \cite{candes2015phase}, a Wirtinger gradient flow algorithm is applied to solve \eqref{1.4}, and it was proven that the algorithm is guaranteed to converge to the global minimizer of \eqref{1.4} when the algorithm is initialized by the so-called spectral initialization and the frame is random Gaussian or Fourier with random masks. To further improve the efficiency and robustness to noise, variants of the Wirtinger flow algorithm are proposed in, e.g., \cite{chen2015solving,zhangs2017median}, and their convergence to the correct solution are provided. More recently, it is revealed  in \cite{sun2016geometric} that the objective in \eqref{1.4} actually has no spurious local minima if $\mathcal{F}$ is a random Gaussian frame with $N\geq O(d\log^3 d)$. Therefore, there are many other efficient algorithms that may be able to find the global minimizer of \eqref{1.4}. However, all these non-convex algorithms assumes random frames, which may be impractical in real applications.

In the fourth category one strives to build special frames with far fewer elements but which still allow for fast and robust reconstruction. In \cite{alexeev2014phase} a deterministic graph-theoretic construction of a frame with $Cd$ measurements was obtained based. This is the few known deterministic construction that uses only $O(d)$ measurements and can robustly reconstruct {\em all} $\vx\in\C^d$, at least in theory. Unfortunately the constant $C$ is very large, so again computationally it would be impractical for large $d$. In \cite{iwen2015robust} a highly efficient phase retrieval scheme using a small number of random linear combinations of Fourier transform measurements is developed. It uses $O(d \log d)$ measurements to guarantee robustness with high probability, and achieves the computational complexity of $O(d\log d)$. In numerical tests it easily performed robust phase retrieval for $d=64,000$ in seconds on a laptop. A drawback is that it is robust only in the probabilistic sense; for a given $\vx$ there is a small probability that the scheme will fail.

In this paper we develop an algorithm for phase retrieval that is both highly efficient and works for very general measurement matrices. 
Our algorithm is based on the ideas of convex relaxation in PhaseLift and
the alternating minimization algorithm used for low rank matrix completion. 
By splitting the variables, our algorithm solves a bi-variant optimization problem whose objective is quadratic in one of the variables with the other fixed.
We shall present both theoretical and numerical results, and discuss its efficient implementation.

\section{Rank One Minimization for Phase Retrieval}
\setcounter{equation}{0}

Let ${\mathcal X} = \{\vx\vx^*: \vx\in\H^d\}$. As we noted in (\ref{1.2}), $\F=\{\vf_n\}$ in $\H^d$ is phase retrievable if and only if the map 
\begin{equation} \label{2.1} 
      \MM_\F(X) = [\vx^*A_1\vx, \dots, \vx^*A_N\vx]^T
                = [\vf_1X\vf_1^*, \dots, \vf_NX\vf_N^*]^T
\end{equation}
is an injective map from ${\mathcal X}$ to $\R^N$, where $A_n = \vf_n\vf_n^*$. In the PhaseLift scheme the phase retrieval problem of solving for $\MM_\F(X) = \vb$ subject to the constraints $X \geq 0$ and has rank 1 (equivalent to $X\in {\mathcal X}$) is being relaxed to solving the convex problem (\ref{1.3})
\begin{equation*} 
   \argmin_X \tr(X) \mhsp\mbox{subject to} \mhsp X\geq 0, ~~\MM_\F(X) = \vb.
\end{equation*}
This relaxation yields the same solution to the original phase retrieval problem with high probability provided that the measurement matrix $A$ is a Gaussian random $N\times d$ matrix with $N=O(d)$ for some unspecified constant, or the DFT matrix with random masks.

Still the drawback is that the measurement matrices are restricted to some specific types, which may or may not be practical for any given application. The optimization requires the use of semi-definite programming, which is slow in general and impractical for phase retrieval for large dimensions. Here we propose a new approach that resolves these difficulties.

The main idea is to relax the requirement $X\in {\mathcal X}$ to simply $\rank(X)=1$. In other words we drop the requirement that $X$ is Hermitian and positive semi-definite. Thus we consider solving the problem
\begin{equation} \label{2.2} 
      \MM_\F(X) = [\vf_1X\vf_1^*, \dots, \vf_NX\vf_N^*]^T=\vb
      \mhsp \mbox{subject to}\mhsp \rank(X)=1,
\end{equation}
or alternatively, given that noise might be present, solving the following problem:
\begin{equation} \label{2.3} 
      \argmin_X \left\|\MM_\F(X)-\vb\right\|
      \mhsp \mbox{subject to}\mhsp \rank(X)=1.
\end{equation}
Observe that in general the solution to (\ref{2.3}) is not unique. If $X$ is a solution then so is $X^*$. To account for this ambiguity we shall use ${\mathcal R}_d(\H)$ to denote the set of $d\times d$ rank one matrices with the equivalence relation $X \equiv X^*$. We shall also let ${\mathcal S}_d(\H)$ denote the set of $d\times d$ Hermitian rank one matrices with entries in $\H$.

\begin{theo}  \label{theo-2.1}
    For $\F= \{\vf_n\}_{n=1}^N \subset \H^d$ and $X\in M_d(\H)$ let 
\begin{equation*} 
      \MM_\F(X) = [\vf_1X\vf_1^*, \dots, \vf_NX\vf_N^*]^T.
\end{equation*}
\begin{itemize}
\item[\rm (A)]~ For a generic $\F \subset \H^d$, $\MM_\F$ is injective on ${\mathcal R}_d(\H)$ if $N \geq 4d-1$ for $\H=\R$, or if $N \geq 8d-3$ for $\H=\C$.

\item[\rm (B)]~ For a generic $\F \subset \H^d$, $\MM_\F$ is injective on ${\mathcal S}_d(\H)$ if $N \geq 2d+1$ for $\H=\R$, or if $N \geq 4d-1$ for $\H=\C$.
\end{itemize}
\end{theo}
\Proof~~
We shall identify $\F=\{\vf_1, \vf_2, \dots, \vf_m\}$ with its frame matrix $F$ whose columns are $\{\vf_n\}$. Consider the set of all 3-tuples
$$
    {\mathcal A}\,:=\, \{(F,X,Y)\}
$$
where $X,Y \in {\mathcal R}_d(\H)$ or $X,Y \in {\mathcal S}_d(\H)$ are distinct and satisfy $\MM_\F(X) = \MM_\F(Y)$. We follow the technique in \cite{balan2006signal} of local dimension counting to prove our theorem.

Let $\H^d_+$ denote the set of vectors of $\H^d$ whose the first nonzero entry is real and positive. Note that any $d\times d$ rank one matrix $Z$ can be written uniquely as $Z = a \vf \vg^*$ where $a\in\H$, $\vf, \vg\in \H^d_+$ and $\|\vf\|=\|\vg\|=1$. Under this factorization $Z\in {\mathcal S}_d(\H)$ if and only if $\vf =\vg$ and $a\in\R$.

To prove the theorem there are 4 cases to be considered, with $\H =\R$ or $\C$ and $X, Y \in {\mathcal R}_d(\H)$ or ${\mathcal S}_d(\H)$. We deal with each case. Due to the similarity of the arguments we shall skip some redundant details.

\vspace{2mm}
\noindent
{\em Case 1:~~$\H=\R$ and  $X,Y \in {\mathcal R}_d(\R)$.}

In this case, because $X, Y$ are distinct in ${\mathcal R}_d(\R)$ each equality $\vf_n X\vf_n^* = \vf_n Y\vf_n^*$ yields a nontrivial constraint in the form of a quadratic polynomial equation for the (real) entries of $\vf_n$. Furthermore, for different $n$ the entries $\vf_n$ are independent variables. Thus viewing the entries of $F$ as points in $\R^{Nd}$, for any distinct $X,Y \in {\mathcal R}_d(\R)$, those satisfying the constraint $\MM_\F(X)=\MM_\F(Y)$ is a real algebraic variety of co-dimension $Nd-N$. By the unique factorization $X =a \vf \vg^*$ discussed above each $X$ has $2d-1$ degrees of freedom. The same $2d-1$ degree of freedom holds also for $Y$. Thus the projection of ${\A}=\{(F,X,Y)\}$ to the first component has local dimension everywhere at most $Nd-N+2(2d-1)$. Suppose that $N \geq 4d-1$. Then this local dimension has
$$
    Nd -N +2(2d-1) \leq Nd-1 <Nd.
$$
In other words, a generic $F\in\R^{N\times d}$ is not a projection of an element in ${\mathcal A}$ to the first component. Thus for a generic $\F$ with $N \geq 4d-1$ the map $\MM_\F$ is injective on ${\mathcal R}_d(\R)$.

\vspace{2mm}
\noindent
{\em Case 2:~~$\H=\R$ and  $X,Y \in {\mathcal S}_d(\R)$.}

All arguments from Case 1 carry to this case, except in the counting of degrees of freedom for $X$ and $Y$. Because now $X=a\vf \vf^*$ there are exactly $d$ degrees of freedom for $X$. The same holds true for $Y$. Thus the projection of ${\mathcal A}=\{(F,X,Y)\}$ to the first component has local dimension everywhere at most $Nd-N+2d$. Suppose that $N \geq 2d+1$. Then this local dimension has
$$
    Nd -N +2d \leq Nd-1 <Nd.
$$
In other words, a generic $F\in\R^{N\times d}$ is not a projection of an element in ${\mathcal A}$ to the first component. Thus for a generic $\F$ with $N \geq 2d+1$ the map $\MM_\F$ is injective on ${\mathcal S}_d(\R)$.

\vspace{2mm}
\noindent
{\em Case 3:~~$\H=\C$ and  $X,Y \in {\mathcal R}_d(\C)$.}

The main arguments from Case 1 carry to this case with slight modifications. A key difference is that we now view $F$ as a point in $\R^{2Nd}$. Each constraint $\vf_n X\vf_n^* = \vf_n Y\vf_n^*$ where $X,Y\in {\mathcal R}_d(\C)$ are distinct now yields an independent nontrivial {\em real} quadratic equation for the real variables ${\rm Re}(\vf_n)$, ${\rm Im}(\vf_n)$. Each $X=a\vf \vg^*$ with $a\in\C$, $\vf,\vg \in \H_+^d$ and $\|\vf\|=\|\vg\|=1$ has 
$2+(2d-2) + (2d-2) = 4d-2$ real degrees of freedom. The same holds for $Y$. Thus the projection of ${\mathcal A}=\{(F,X,Y)\}$ to the first component has real local dimension everywhere at most $2Nd-N+2(4d-2)$. Suppose that $N \geq 8d-3$. Then this real local dimension has
$$
    2Nd -N +2(4d-2) \leq 2Nd-1 <2Nd.
$$
In other words, a generic $F\in\C^{N\times d}$ is not a projection of an element in ${\mathcal A}$ to the first component. Thus for a generic $\F$ with $N \geq 8d-3$ the map $\MM_\F$ is injective on ${\mathcal R}_d(\C)$.

\vspace{2mm}
\noindent
{\em Case 4:~~$\H=\C$ and  $X,Y \in {\mathcal S}_d(\C)$.}

All arguments from Case 3 carry to this case, except in the counting of degrees of freedom for $X$ and $Y$. Because now $X=a\vf \vf^*$ where $a\in\R$, $\vf\in \H_+^d$ and $\|\vf\|=1$ there are exactly $1+2d-2=2d-1$ real degrees of freedom for $X$. The same holds true for $Y$. Thus the projection of ${\mathcal A}=\{(F,X,Y)\}$ to the first component has real local dimension everywhere at most $2Nd-N+2(2d-1)$. Suppose that $N \geq 4d-1$. Then this local dimension has
$$
    Nd -N +2(2d-1) \leq Nd-1 <Nd.
$$
In other words, a generic $F\in\C^{N\times d}$ is not a projection of an element in ${\mathcal A}$ to the first component. Thus for a generic $\F$ with $N \geq 4d-1$ the map $\MM_\F$ is injective on ${\mathcal S}_d(\C)$.
\eproof

We can now reformulate the phase retrieval problem into two alternative optimization problems. Each rank one matrix $X$ can be written as $X=\vx \vy^*$ for some $\vx,\vy \in\H^d$, although this representation is not unique. We have

\begin{theo}  \label{theo-2.2}
     Let $\F=\{\vf_n\}_{n=1}^N$ be vectors in $\H^d$ such that $\MM_\F$ is injective on ${\mathcal R}_d(\H)$. Let $\vx_0\in\H^d(\H)$ and $\vb=\MM_\F(\vx_0\vx_0^*) =[|\inner{\vf_1,\vx_0}|^2, \dots, |\inner{\vf_N,\vx_0}|^2]^T$. Then any global minimizer
\begin{equation}  \label{2.4}
      (\hat\vx, \hat\vy) = \argmin_{\vx,\vy\in\H^d} \left\|\MM_\F(\vx\vy^*)-\vb\right\|
\end{equation}
must satisfy $\hat\vx=\hat\vy = c\vx_0$ for some $|c|=1$. 
\end{theo}

\Proof   The result follows trivially from the injectivity of $\MM_\F(X)$ on ${\mathcal R}_d(\H)$. Clearly if $\hat\vx=\hat\vy = c\vx_0$ then $\vx\vy^*=\vx_0\vx_0^*$ which gives the global minimizer. Conversely, the global minimizer must have $\MM_\F(\vx\vy^*)=\vb$. The injectivity now implies that $\hat\vx\hat\vy = \vx_0\vx_0^*$. Thus $\hat\vx=\hat\vy = c\vx_0$ for some $|c|=1$. 
\eproof

The above minimization problem is not convex so solving for the global minimum is very challenging. Such is the case for solving the phase retrieval problem in general. The advantage of the above formulation is that it allows us to use the popular alternating minimization technique used for many other applications such as low rank matrix completion, see e.g. \cite{jain2013low,hardt2014understanding,wei2016guarantees,tanner2016low} the references therein. In the alternating minimization algorithm, we first pick an initial $\vx_1$ and minimize $\|\MM_\F(\vx_1\vy^*)-\vb\|$ with respect to $\vy$ to obtain $\vy_1$. This step is a standard $\ell_2$-minimization and is linear problem. From $\vy_1$ we then update $\vx$ to $\vx_2$ via minimizing $\|\MM_\F(\vx\vy_1^*)-\vb\|$. This process is iterated to yield a sequence $\vx_k\vy_k^*$. Often the sequence converges to the desired result.

The drawback of the above setup is that because there is no penalty for $\vx\vy^*$ being non-symmetric, when noise is added to the measurement vector $\vb$, the stability and robustness is harder to analyze. It also requires, at least in theory, almost twice as many measurements as the minimally required number for phase retrieval. A better alternative is to add a regularization term to the previous minimization problem. Let $\lambda>0$ and 
\begin{equation}  \label{2.5}
      E_{\lambda,\vb}(\vx,\vy)= \left\|\MM_\F(\vx\vy^*)-\vb\right\|^2 +\lambda \|\vx-\vy\|^2.
\end{equation}
Below we study the consequences of minimizing this function. In particular, we wish to establish certain robustness properties.

\begin{lem} \label{lem-upperbound}
     Let $\F=\{\vf_n\}_{n=1}^N$ be a frame in $\H^d$. Let $X\in M_d(\H)$. Then $\|\MM_\F(X)\|_1 \leq C\|X\|_*$ where $\|X\|_*$ denotes the nuclear norm of $X$, and $C$ is the upper frame bound of $\F$, i.e. $C$ is the largest eigenvalue of $FF^*$ where $F= [\vf_1, \dots, \vf_N]$ is the frame matrix for $\F$. Furthermore, this $C$ is optimal.
\end{lem}

\Proof  Assume that $X=\vv\vv^*$ for some $\vv\in\H^d$ and $\|\vv\|=1$. Then $\|X\|_*=1$
$$
    \|\MM_\F(X)\|_1 = \sum_{n=1}^N |\inner{\vf_n,\vv}|^2 \leq C\|\vv\|^2 = C.
$$
Note that here this constant $C$ is the best possible since it can be achieved by taking $\vv$ to be an eigenvector of $FF^*$ corresponding to its largest eigenvalue. Now assume that $X$ is Hermitian. Then we may write $X$ as
$$
     X = \sum_{j=1}^d \lambda_j \vv_j\vv_j^*
$$
where $\{\vv_j\}$ is an orthonormal basis for $\H^d$. Thus
$$
     \|\MM_\F(X)\|_1 \leq \sum_{j=1}^d |\lambda_j|\|\MM_\F(\vv_j\vv_j^*)\|_1
     \leq C \sum_{j=1}^d |\lambda_j| = C\|X\|_*.
$$
For a non-Hermitian $X$, let $Y= \frac{1}{2}(X+X^*)$. Then $\|Y\|_* \leq \|X\|_*$, and
$$
    \|\MM_\F(X)\|_1 = \|\MM_\F(Y)\|_1 \leq C\|Y\|_* \leq C \|X\|_*.
$$
\eproof

For a phase retrievable set $\F=\{\vf_n\}_{n=1}^N$ in $\H^d$ we say $\MM_\F$ satisfies the {\em $c$-stability condition} if for any $\vx, \vy\in\H^d$ we have
\begin{equation}  \label{2.6}
  \left\|\MM_\F(\vx\vx^*)-\MM_\F(\vy\vy^*)\right\| 
         \geq c \|\vx\vx^* -\vy\vy^*\|_*.
\end{equation}

\begin{theo}  \label{theo-2.3}
     Let $\F=\{\vf_n\}_{n=1}^N$ be phase retrievable in $\H^d$. Let $\vx_0\in\H^d(\H)$ and $\vb=\MM_\F(\vx_0\vx_0^*) =[|\inner{\vf_1,\vx_0}|^2, \dots, |\inner{\vf_N,\vx_0}|^2]^T$. Then 
\begin{itemize} 
\item[\rm (A)] Any global minimizer
$$
      (\hat\vx, \hat\vy) = \argmin_{\vx,\vy\in\H^d} E_{\lambda,\vb}(\vx,\vy)
$$
must satisfy $\hat\vx \hat\vy^*=\vx_0\vx_0^*$, or equivalently $\hat\vx=\hat\vy = c\vx_0$ for some $|c|=1$.
\item[\rm (B)] Let $\vb'\in\H^d$ such that $\|\vb-\vb'\| \leq \ep$. Assume that $\MM_\F$ satisfies the $c$-stability condition for some $c>0$. Then any $\vx, \vy\in\H^d$ such that $E_{\lambda,\vb'}(\vx,\vy)\leq \delta^2$ must satisfy 
\begin{equation}  \label{2.7}
    \left\|\vz \vz^* -\vx_0\vx_0^*\right\|_* \leq 
         \frac{1}{c}\Bigl(\frac{C}{4\lambda}\delta^2 + \delta+\ep\Bigr), 
\end{equation}
where $\vz=\frac{1}{2}(\vx+\vy)$ and $C$ is the upper frame bound of $\F$. 
\end{itemize}
\end{theo}
\Proof  Part (A) is rather straightforward. Note that $E_{\lambda,\vb}(\vx_0,\vx_0)=0$, so we must have $E_{\lambda,\vb}(\hat\vx,\hat\vy)=0$. It follows that $\hat\vx=\hat\vy$. Hence $\MM_\F(\hat\vx\hat\vx^*)=\vb =\MM_\F(\vx_0\vx_0^*)$. The fact that $\F$ is phase retrievable now implies $\hat\vx\hat\vx^*=\vx_0\vx_0^*$.

To prove part (B), we have $\lambda\|\vx-\vy\|^2 \leq \delta^2$. Thus $\|\vx-\vy\| \leq \delta/\sqrt{\lambda}$. Let $Z=\frac{1}{2}(\vx\vy^*+\vy\vx^*)$. Clearly $\MM_\F(Z)=\MM_\F(\vx\vy^*)$. Furthermore one checks easily that
$$
   \vz\vz^* -Z= \frac{1}{4}(\vx-\vy)(\vx-\vy)^*.
$$
Hence $\|\vz\vz^* -Z\|_* =\frac{1}{4}\|\vx-\vy\|^2\leq \frac{\delta^2}{4\lambda}$. It follows that
\begin{align*}
   \|\MM_\F(\vz\vz^*)&-\MM_\F(\vx_0\vx_0^*)\| = \|\MM_\F(\vx\vx^*)-\vb\|  \\
      & \leq \|\MM_\F(\vz\vz^*)-\MM_\F(Z)\| + \|\MM_\F(Z)-\vb'\|+\|\vb'-\vb\| \\
      & \leq \|\MM_\F(\vz\vz^*)-\MM_\F(Z)\|_1 + \|\MM_\F(Z)-\vb'\|+\|\vb'-\vb\| \\
       & \leq C\frac{\delta^2}{4\lambda} + \delta+\ep.
\end{align*}
The $c$-stability condition now implies (\ref{2.7}) immediately.
\eproof

\section{Alternating Minimization Algorithm}
\setcounter{equation}{0}
From the formulation in the previous section, the phase retrieval problem is solved robustly by finding a global minimizer of  $E_{\lambda,\vb}$ in \eqref{2.5}. This section is devoted to fast algorithms for solving such a minimization problem. In Section \ref{algorithm}, we introduce a fast alternating gradient descent algorithm for $\min_{\vx,\vy}E_{\lambda,\vb}(\vx,\vy)$. In Section \ref{convergence}, we prove the convergence of the proposed algorithm. 

\subsection{Alternating gradient descent algorithm}
\label{algorithm}
Since $\lambda$ and $\vb$ are fixed during the minimization procedure, we drop the subscripts in $E_{\lambda,\vb}$ for simplicity. That is, we solve
\begin{equation}
\label{model_now}
\underset{\vx,\vy}{\arg\min}\ E(\vx,\vy),
\end{equation}
where
\begin{equation}
\label{E_function}
E(\vx,\vy)=\frac{1}{m}\sum^N_{n=1}|\vx^*\vf_n\vf^*_n\vy-b_n|^2+\lambda\|\vx-\vy\|^2.
\end{equation}
Since the first term in $E(\vx,\vy)$ is quartic in $(\vx,\vy)$, Eq. \eqref{model_now} is a non-convex optimization. However, when one of the variables $\vx$ or $\vy$ is fixed, $E(\vx,\vy)$ is quadratic with respect to the other variable. Therefore, it is natural to solve \eqref{model_now} by an alternating scheme.

We use the following alternating gradient descent algorithm: Fixing $\vx$, we minimize $E(\vx,\vy)$ with respect to $\vy$ by one step of gradient descent, and vice versa. More precisely, we define, for $k=0,1,2,\cdots$,
\begin{equation}
\label{iteration_scheme}
\left\{\begin{array}{l}
\vx_{k+1}=\vx_{k}-\alpha_k\nabla_{\vx} E(\vx_k,\vy_k),\\
\vy_{k+1}=\vy_k-\beta_k\nabla_{\vy} E(\vx_{k+1},\vy_k),
\end{array}
\right.
\end{equation}
where $\alpha_k$ and $\beta_k$ are step sizes. Since $E$ is a real-valued function with complex variables, the gradients $\nabla_{\vx} E$ and $\nabla_{\vy} E$ in \eqref{iteration_scheme} are in the sense of Wirtinger gradient \cite{candes2015phase}.

Since the gradient descent is applied to only one of the variables $\vx$ and $\vy$, the corresponding Hessian matrix have a much smaller norm than the Hessian of $E$ with respect to $(\vx,\vy)$. Consequently, a much larger step size is allowed in the alternating gradient descent than the standard gradient descent for minimizing \eqref{E_function}, which leads to a faster convergence. The alternating gradient descent algorithm is also faster than the Wirtinger flow (WF) \cite{candes2015phase} algorithm, where $G(\vx)=\frac{1}{N}\sum^N_{n=1}(|\vf^*_n\vx|^2-b_n)^2$ is minimized via a gradient flow. As explained in Appendix A, in the real case, when the iterates $\vx$ and $\vy$ are sufficiently close, our proposed alternating gradient descent algorithm is 1.5 times faster than the WF algorithm in terms of the decreasing of the objectives.


 

The initialization of our proposed algorithm is obtained via a spectral method, which is the same as that in the Wirtinger flow algorithm \cite{candes2015phase}. When $\vf_n$, $n=1,\ldots,N$, follow certain probability distributions (e.g., Gaussian), the expectation of $Y=\frac{1}{N}\sum^N_{n=1}b_n\vf_n\vf^*_n$ has a leading principle eigenvector $\tilde{\vx}$. Therefore, we choose $\vx_0=\vy_0=\vz_0$, where $\vz_0$ is the leading principle eigenvector of $Y$. For completeness, Algorithm \ref{initialization} lists how to calculate the initial guess for our proposed algorithm.

\begin{algorithm}[H]
\KwIn{Observations $b_n,n=1,\cdots, N$.}
\KwOut{Initial guess $\vx_0=\vy_0=\vz_0$.}

Set
$$\theta^2=d\frac{\sum_nb_n}{\sum_n\|\vf_n\|^2},$$ where $\vf_n\in\mathbb{C}^d,n=1,\cdots,N$ is the sampling vectors\;
$\vz_0$ is the eigenvector corresponding to the largest eigenvalue of 
$$Y=\frac{1}{N}\sum^N_{n=1}b_n\vf_n\vf^*_n$$ and $\|\vz_0\|=\theta$.
\caption{Initialization.}
\label{initialization}
\end{algorithm}

\subsection{Convergence}
\label{convergence}
In this section, we will show the convergence of the alternating gradient descent algorithm \eqref{iteration_scheme}. 
More precisely, for any initial guess, we prove that algorithm \eqref{iteration_scheme} converges to a critical point of $E$.  

We first present a lemma, which shows the coercivity of $E$.
\begin{lem}\label{lem_coercive}
If $F$ is of full rank, i.e., $\mathrm{rank}(F)=d$, then the function $E(\vx,\vy)$ is coercive, i.e., $E(\vx,\vy)\to\infty$ as $\|(\vx,\vy)\|\to\infty$.
\end{lem}

\Proof
Since $F$ is of full row rank, there exists a constant $C_1$ such that, for any $\vx$,
$$
\||F^*\vx|^2\|=\|F^*\vx\|_4^2\geq C_0\|F^*\vx\|^2\geq C_1\|\vx\|^2.
$$
Also, there exists a constant $C_2$ such that, for any $\vx$ and $\vz$,
$$
\|(\overline{F^*\vx})\circ(F^*\vz)\|\leq\|F^*\vx\|\|F^*\vz\|\leq C_2\|\vx\|\|\vz\|,
$$
where $\circ$ is the componentwise product.

Let $\|(\vx,\vy)\|= M$. If $\|\vx\|\leq \frac{M}{2}$, then $\|\vy\|\geq\frac{\sqrt3 M}{2}$ and$$E(\vx,\vy)\geq\lambda\|\vx-\vy\|^2\geq \frac{\lambda(\sqrt3-1)^2M}{4}.$$ Similarly, if $\|\vy\|_2\leq\frac{M}{2}$, then $E(\vx,\vy)\geq\frac{\lambda(\sqrt3-1)^2M}{4}.$
Otherwise, both $\|\vx\|_2> \frac{M}{2}$ and $\|\vy\|_2>\frac{M}{2}$.
Define $\vz=\vy-\vx$. In this case, if $\|\vz\|\leq \frac{C_1M}{4C_2}$, then
\begin{equation*}
\begin{split}
E(\vx,\vy)&\geq\frac{1}{N}\sum^N_{n=1}|\vx^*\vf_n\vf^*_n\vy-b_n|^2=
\frac{1}{N}\|(\overline{F^*\vx})\circ(F^*(\vx+\vz))-\vb\|^2\cr
&\geq\frac{1}{N}(\||F^*\vx|^2\|-\|\vb\|-\|(\overline{F^*\vx})\circ(F^*\vz)\|)^2\geq\frac{1}{N}(C_1\|\vx\|^2-\|\vb\|-C_2\|\vx\|\|\vz\|)^2\cr
&=\frac1N(\|\vx\|(C_1\|\vx\|-C_2\|\vz\|)-\|\vb\|)^2\geq \frac1N(M^2/8-\|\vb\|)^2,
\end{split}
\end{equation*}
otherwise $E(\vx,\vy)\geq\lambda\|\vz\|^2\geq \frac{\lambda C_1^2M^2}{16C_2^2}$. In all the cases, as $M\to\infty$, the lower bounds approach infinity.
\eproof

Now we can prove the convergence of \eqref{iteration_scheme}.
\begin{theo}
Assume $\rank(F)=d$. Then, for any initial guess $(\vx_0,\vy_0)$, the sequence $\{(\vx_k,\vy_k)\}_k$ generated by \eqref{iteration_scheme} with a suitable step size
converges to a critical point of $E$.
\end{theo}

\Proof
For simplicity, we assume $\alpha_k=\beta_k=\gamma$ for all $k$. The proof with variant step sizes can be done similarly. 
Since $E(\vx,\vy)$ is a quadratic function with respect to $\vx$,
Taylor's expansion gives
\begin{equation}\label{eq:objdecay1}
\begin{split}
E(&\vx_{k+1},\vy_k)=E(\vx_k,\vy_k)+\begin{bmatrix}\vx_{k+1}-\vx_k\\\overline{\vx_{k+1}}-\overline{\vx_k}\end{bmatrix}^*\begin{bmatrix}\nabla_{\vx}E(\vx_k,\vy_k)\\\nabla_{\overline{\vx}}E(\vx_k,\vy_k)\end{bmatrix}\\
&\qquad\qquad\qquad+\frac{1}{2}\begin{bmatrix}\vx_{k+1}-\vx_k\\\overline{\vx_{k+1}}-\overline{\vx_k}\end{bmatrix}^*\begin{bmatrix}\nabla^2_{\vx\vx}E(\vx_k,\vy_k)&\nabla^2_{\vx\overline{\vx}}E(\vx_k,\vy_k)\\
\nabla^2_{\overline{\vx}\vx}E(\vx_k,\vy_k)&\nabla^2_{\overline{\vx}\overline{\vx}}E(\vx_k,\vy_k)\end{bmatrix}\begin{bmatrix}\vx_{k+1}-\vx_k\\\overline{\vx_{k+1}}-\overline{\vx_k}\end{bmatrix}\\
&=E(\vx_k,\vy_k)-2\gamma\|\nabla_{\vx}E(\vx_k,\vy_k)\|^2+\gamma^2\nabla_{\vx}E(\vx_k,\vy_k)^*\nabla^2_{\vx\vx}E(\vx_k,\vy_k)\nabla_{\vx}E(\vx_k,\vy_k)\\
&=E(\vx_k,\vy_k)-2\gamma\left(\|\nabla_{\vx}E(\vx_k,\vy_k)\|^2-\frac{\gamma}{2}\nabla_{\vx}E(\vx_k,\vy_k)^*\nabla^2_{\vx\vx}E(\vx_k,\vy_k)\nabla_{\vx}E(\vx_k,\vy_k)\right).
\end{split}
\end{equation}
Similarly, because $E(\vx,\vy)$ is a quadratic function with respect to $\vy$,
by Taylor's expansion, we obtain
\begin{equation}
\label{eq:objdecay2}
\begin{split}
&E(\vx_{k+1},\vy_{k+1})
=E(\vx_{k+1},\vy_k)-\cr&\qquad2\gamma\left(\|\nabla_{\vy}E(\vx_{k+1},\vy_k)\|^2-\frac{\gamma}{2}\nabla_{\vy}E(\vx_{k+1},\vy_k)^*\nabla^2_{\vy\vy}E(\vx_{k+1},\vy_k)\nabla_{\vy}E(\vx_{k+1},\vy_k)\right).
\end{split}
\end{equation}
By Lemma \ref{lem_coercive}, the level set $\mathcal{S}=\{(\vx,\vy):E(\vx,\vy)\leq E(\vx_0,\vy_0)\}$ is a bounded closed set. Therefore, the continuous functions $\nabla^2_{\vx\vx}E(\vx,\vy)$ and $\nabla^2_{\vy\vy}E(\vx,\vy)$ are bounded on $\mathcal{S}$. Let $M>0$ be the bound, i.e.,
\begin{eqnarray*}
\|\nabla^2_{\vx\vx}E(\vx,\vy)\|\le M, \quad \|\nabla^2_{\vy\vy}E(\vx,\vy)\|\le M,\quad\forall\ (\vx,\vy)\in\mathcal{S}.
\end{eqnarray*}
Suppose $(\vx_k,\vy_k)\in\mathcal{S}$. 
Choose $\gamma\in(0,2/M)$,
so that \eqref{eq:objdecay1} and \eqref{eq:objdecay2} implies that $E(\vx_{k+1},\vy_{k+1})\leq E(\vx_k,\vy_k)$ and
\begin{eqnarray}
\label{lower_bound}
E(\vx_k,\vy_k)-E(\vx_{k+1},\vy_{k+1})&=&E(\vx_k,\vy_k)-E(\vx_{k+1},\vy_k)+E(\vx_{k+1},\vy_k)-E(\vx_{k+1},\vy_{k+1})\nonumber\\
&\ge&\zeta(\|\nabla_{\vx}E(\vx_k,\vy_k)\|^2+\|\nabla_{\vy}E(\vx_{k+1},\vy_{k})\|^2)\\
&=&\frac{\zeta}{\gamma}(\|\vx_{k+1}-\vx_k\|^2+\|\vy_{k+1}-\vy_k\|^2)\nonumber
\end{eqnarray}
with $\zeta=2\gamma\left(1-\frac{\gamma M}{2}\right)$ .
Therefore, $(\vx_{k+1},\vy_{k+1})\in\mathcal{S}$. Thus, by induction, $(\vx_{k},\vy_{k})\in\mathcal{S}$ and \eqref{lower_bound} hold for all $k$ as long as $\gamma\in(0,2/M)$. 

Summing \eqref{lower_bound} over $k$ from $0$ to $+\infty$, we obtain
$$
E(\vx_{0},\vy_{0})-\lim_{k\to+\infty}E(\vx_{k},\vy_{k})\geq
\zeta\sum_{k=0}^{+\infty}(\|\nabla_{\vx}E(\vx_k,\vy_k)\|^2+\|\nabla_{\vy}E(\vx_{k+1},\vy_{k})\|^2).
$$
Because $E(\vx_{k},\vy_{k})\geq 0$ is monotonically nonincreasing according to \eqref{lower_bound}, its limit exists and is finite, which implies
\begin{eqnarray*}
&\lim_{k\to +\infty}\|\nabla_{\vx}E(\vx_k,\vy_k)\|=\lim_{k\to +\infty}\|\vx_{k+1}-\vx_k\|=0,\\
&\lim_{k\to +\infty}\|\nabla_{\vy}E(\vx_{k+1},\vy_k)\|
=\lim_{k\to +\infty}\|\vy_{k+1}-\vy_k\|=0.
\end{eqnarray*}
This, together with the continuity of $\nabla_{\vx}E$, $\nabla_{\vy}E$, and the norm function, means that any clustering point of $\{(\vx_k,\vy_k)\}_k$ is a critical point of $E(\vx,\vy)$.

It remains to prove that $\{(\vx_k,\vy_k)\}_k$ is convergent, which is done by checking that $\{(\vx_k,\vy_k)\}_k$ is a Cauchy sequence.
Since $E(\vx,\vy)$ is a real-valued polynomial function, it belongs to a semi-algebraic set. By \cite[Theorem 3]{bolte2014proximal}, there exists a differentiable and concave function $\psi(t)$ such that
\begin{equation}
\label{kl}\psi'(E(\vx_k,\vy_k)-E(\hat{\vx},\hat{\vy}))\cdot\left\|\begin{bmatrix}\nabla_{\vx}E(\vx_k,\vy_k)\\\nabla_{\vy}E(\vx_k,\vy_k)\end{bmatrix}\right\|\ge1
\end{equation}
for any $k$ and for any critical point $(\hat{\vx},\hat{\vy})$ of $E$.
Since $\psi(t)$ is concave, by the inequalities \eqref{lower_bound} and \eqref{kl}, we have 
\begin{equation}\label{kl2}
\begin{split}
\psi(E(\vx_k,\vy_k)&-E(\hat{\vx},\hat{\vy}))-\psi(E(\vx_{k+1},\vy_{k+1})-E(\hat{\vx},\hat{\vy}))\\
\ge&\psi'(E(\vx_k,\vy_k)-E(\hat{\vx},\hat{\vy}))(E(\vx_k,\vy_k)-E(\vx_{k+1},\vy_{k+1}))\\
\ge&\frac{\zeta}{\gamma}\frac{\|\vx_{k+1}-\vx_k\|^2+\|\vy_{k+1}-\vy_k\|^2}{\sqrt{\|\nabla_{\vx}E(\vx_k,\vy_k)\|^2+\|\nabla_{\vy}E(\vx_k,\vy_k)\|^2}}.
\end{split}
\end{equation}
Furthermore, 
\begin{equation}
\label{inequality1}
\|\nabla_{\vx}E(\vx_{k},\vy_{k})\|
=\frac{1}{\gamma}\left\|(\vx_{k+1}-\vx_{k})\right\|
\end{equation}
and 
\begin{equation}
\label{inequality2}
\begin{split}
\|\nabla_{\vy}E(\vx_{k},\vy_{k})\|
&=\left\|\nabla_{\vy}E(\vx_{k},\vy_{k})-\nabla_{\vy}E(\vx_{k+1},\vy_{k})+\frac{1}{\gamma}(\vy_{k+1}-\vy_{k})\right\|\\
&\le\left\|\nabla_{\vy}E(\vx_{k},\vy_{k})-\nabla_{\vy}E(\vx_{k+1},\vy_{k})\right\|+\frac{1}{\gamma}\left\|\vy_{k+1}-\vy_{k}\right\|\\
&\le M'\left\|\vx_{k+1}-\vx_{k}\right\|+\frac{1}{\gamma}\left\|\vy_{k+1}-\vy_k\right\|,
\end{split}
\end{equation}
where $M'=\sup_{(\vx,\vy)\in\mathcal{S}}\|\nabla^2_{\vx\vy}E(\vx,\vy)\|$ that is finite. Plugging \eqref{inequality1} and \eqref{inequality2} into \eqref{kl2} gives
\begin{equation*}\label{kl3}
\begin{split}
\psi(E(\vx_k,\vy_k)&-E(\hat{\vx},\hat{\vy}))-\psi(E(\vx_{k+1},\vy_{k+1})-E(\hat{\vx},\hat{\vy}))\\
\ge&C\left(\|\vx_{k+1}-\vx_k\|^2+\|\vy_{k+1}-\vy_k\|^2\right)^{1/2},
\end{split}
\end{equation*}
where $C=\frac{\zeta}{\gamma M'+1}$.
Summing it over $k$, we get
\begin{equation*}
\label{finite}
\begin{split}
\sum^{+\infty}_{k=0}&\left(\|\vx_{k+1}-\vx_k\|^2+\|\vy_{k+1}-\vy_k\|^2\right)^{1/2}\cr
&\leq \frac{1}{C}
\left(\psi(E(\vx_0,\vy_0)-E(\hat{\vx},\hat{\vy}))-\lim_{k\to\infty}\psi(E(\vx_{k},\vy_{k})-E(\hat{\vx},\hat{\vy}))\right).
\end{split}
\end{equation*}
The right hand side is finite, as $\psi$ is smooth and $\lim_{k\to\infty}E(\vx_k,\vy_k)$ is finite. This verifies that $\{(\vx_k,\vy_k)\}_k$ is a Cauchy sequence, 
and therefore it is convergent.
\eproof

\section{Numerical Implementation}
In this section, we present some numerical experiments to evaluate the proposed alternating gradient descent algorithm and compare it with the Wirtinger Flow (WF) algorithm \cite{candes2015phase}. As demonstrated in Section \ref{synthetic} on synthetic data and Section \ref{real} on real image data, our proposed algorithm is 
more efficient in terms that a smaller number of iterations are required to achieve the same recovery accuracy.

\subsection{Experiment Setup}
All experiments are carried out on a PC with a 3.20 GHz Intel Core i5 Processor and 8GB memory. The initialization of our proposed algorithm is described in Algorithm \ref{initialization}, which is run by 50 iterations of the power method. 
In WF algorithm \cite{candes2015phase}, the step size is chosen heuristically and expementally as $\tilde{\mu}_{\tau}=\min(1-e^{-\tau/\tilde{\tau}_0},\tilde{\mu}_{max})$, with $\tilde{\tau}_0=330$ and $\tilde{\mu}_{max}=0.2\ \mathrm{or}\ 0.4$, 
which is the most efficient according to our test.
Following this, the step size of our method  is also chosen in the form as $\mu_{\tau}=\min(1-e^{-\tau/\tau_0},\mu_{max})$ and the tuning parameter $\lambda_{\tau}=\lambda_0 e^{-\xi\tau}$. 
The parameters $\tilde{\tau}_0$, $\tau_0$, $\tilde{\mu}_{max}$ $\mu_{max}$, $\lambda_0$ and $\xi$ will be specified later.

Throughout the test, we mainly focus on the Gaussian model and the coded diffraction (CDF) model. In  the Gaussian model, we collect the data $b_n=|\vf^*_n\vx|^2$ with the sampling vectors distributed as Gaussian model, that is,
$$
\vf_n\overset{\mathrm{i.i.d.}}{\thicksim}\left\{
\begin{array}{ll}
\mathcal{N}(0,I/2)+i\mathcal{N}(0,I/2), & \mathrm{if}\ \vf_n\in\mathbb{C}^d,\\
\mathcal{N}(0,I), & \mathrm{if}\ \vf_n\in\mathbb{R}^d,\end{array}
\right.
$$
where $\mathcal{N}(0,V)$ is the real mean-zero Gaussian distribution with covariance matrix $V$. 
In the CDF model, we acquire the data via
$$
b_{p,q}=\sum_{p,q}|\vx^*\va_{p,q}|^2,\quad \mathrm{with}\ 0\le q\le d-1,\ 1\le p\le L,
$$
with $\va_{p,q}=G_p\vf_q$, where $\vf^*_q$ is the $q$-th row of the $d\times d$ discrete Fourier Transform (DFT) 
matrix and $G_p$ is a diagonal matrix with i.i.d. diagonal entries 
$g_p(0),~g_p(1),~\cdots$, $g_p(d-1)$ randomly drawn from $\left\{\pm\frac{\sqrt{2}}{2},~\pm\frac{\sqrt{2}}{2}i\right\}$ with probability $\frac15$ for each element, and
$\left\{\pm\sqrt{3},~\pm\sqrt{3}i\right\}$ with probability $\frac1{20}$ for each element.

  \subsection{Synthetic data}
 \label{synthetic}
 In this subsection, we test the algorithms on synthetic data. Following \cite{candes2015phase}, we are interested in the two signals described below:
\begin{itemize}
\item \emph{Random low-pass signals.} The true signal $\tilde{\vx}\in\mathbb{C}^d$ is generated by 
$$
\tilde{x}[t]=\sum^{M/2}_{k=-(M/2-1)}(r_k+ij_k)e^{2\pi i(k-1)(t-1)/d}
$$
where $M=\frac{d}{8}$, and $r_k$ and $j_k$ are i.i.d. obeying the standard normal distribution.
\item \emph{Random Gaussian signals.} The true signal $\tilde{\vx}\in\mathbb{C}^d$ is a random complex Gaussian vector with i.i.d. entries of the form 
$$\tilde{x}[t]=\sum^{d/2}_{k=-(d/2-1)}(r_k+ij_k)e^{2\pi i(k-1)(t-1)/d},$$
where $r_k$ and $j_k$ are i.i.d. normal distribution $\mathcal{N}(0,\frac{1}{8})$.
\end{itemize}

\begin{figure}[h]
\begin{center}
\subfigure[Gaussian model]{
\label{Fig.sub11}
\includegraphics[width=0.45\textwidth]{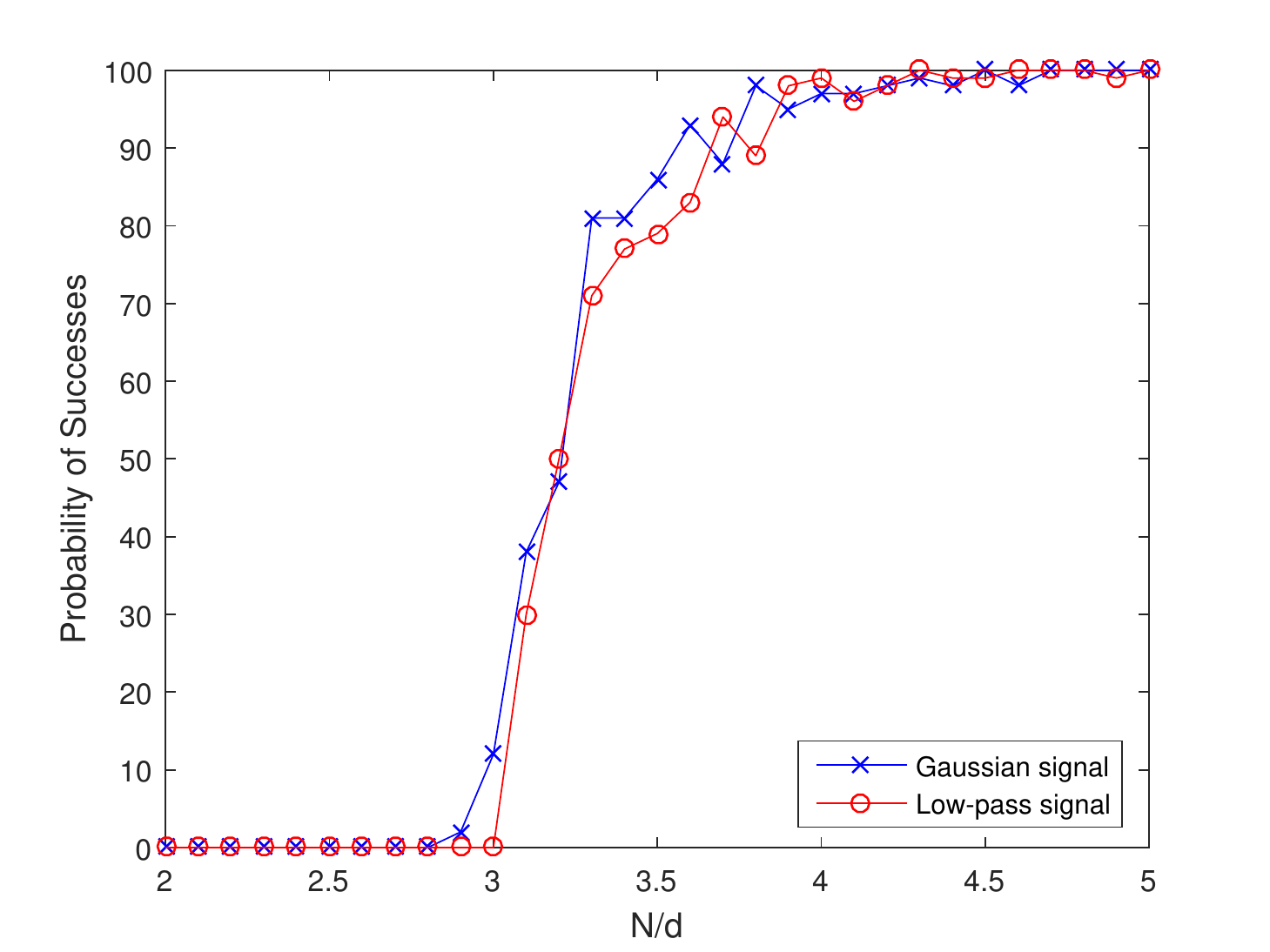}}
\subfigure[Coded diffraction model]{
\label{Fig.sub12}
\includegraphics[width=0.45\textwidth]{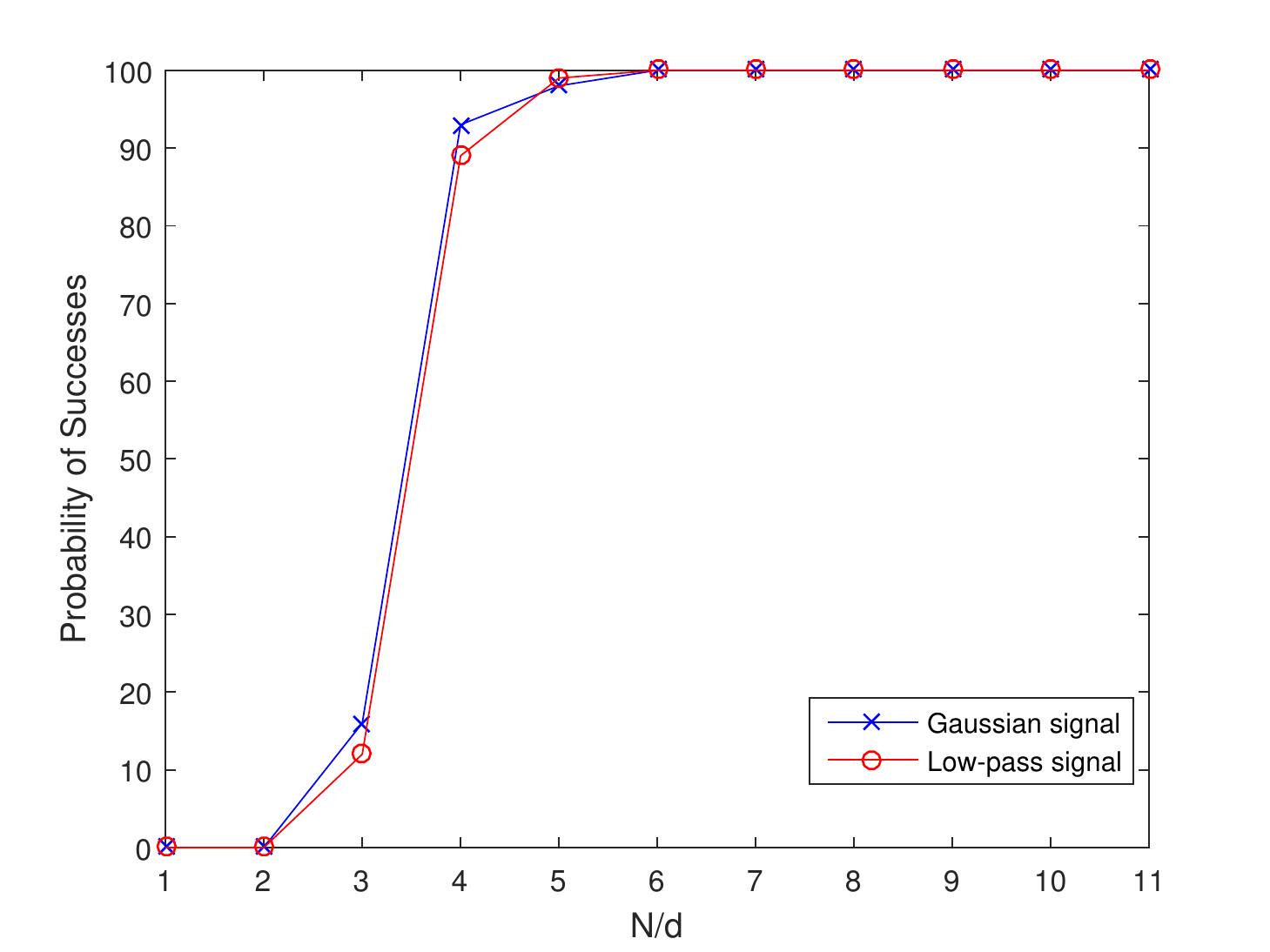}}
\caption{The plot of the probability of success versus $N/d$.}
\label{fig:measurements}
\end{center}
\end{figure}

We first evaluate the effectiveness of our proposed algorithm in terms of the smallest $N$ required for successful phase retrieval. 
We use 100 trials for both the Gaussian and CDF models. In each trial, we generate the random sampling vectors according to the Gaussian or CDF model and stop the alternating iteration after 2500 iterations (1250 iterations for $\vx$ and 1250 iterations for $\vy$ corresponding to our method). We declare it is successful if the relative error of the construction $\mathrm{dist}(\tilde{\vx},\hat{\vx})/\|\tilde{\vx}\|<10^{-5}$, where $\hat{\vx}$ is the numerical solution by our alternating minimization algorithm. The empirical probability of success is defined as the average of success over 100 trials.   
We use $d=128$. In the Gaussian model, we choose $\tau_0=\tilde{\tau}_0=330$, $\tilde{\mu}_{max}=0.2$, $\mu_{max}=0.4$, $\lambda_0=300$ and $\xi=0.15/330$ for random Gaussian signal, and $\tau_0=\tilde{\tau}_0=330$, $\tilde{\mu}_{max}=0.2$, $\mu_{max}=0.4$, $\lambda_0=5$ and $\xi=0.05/300$ for the random low-pass signal. In the CDF model, we choose $\tau_0=\tilde{\tau}_0=330$, $\tilde{\mu}_{max}=0.2$, $\mu_{max}=0.4$, $\lambda_0=0.2$ and $\xi=0.0015/330$ for random Gaussian signal, and $\tau_0=\tilde{\tau}_0=330$, $\tilde{\mu}_{max}=0.2$, $\mu_{max}=0.4$, $\lambda_0=0.05$ and $\xi=1.5/330$ for the random low-pass signal.. We plot the empirical probability of success against the over sampling ratio $N/d$ in Figure \ref{fig:measurements}. We see that the minimum oversampling ratios for an almost $100\%$ successful phase retrieval by our algorithm are around $4.3$ for the Gaussian model and $6$ for the CDF model, which is slightly better or the same as the requirement of the WF algorithm as reported in \cite{candes2015phase}.



\begin{figure}[h]
\begin{center}
\includegraphics[width=0.5\textwidth]{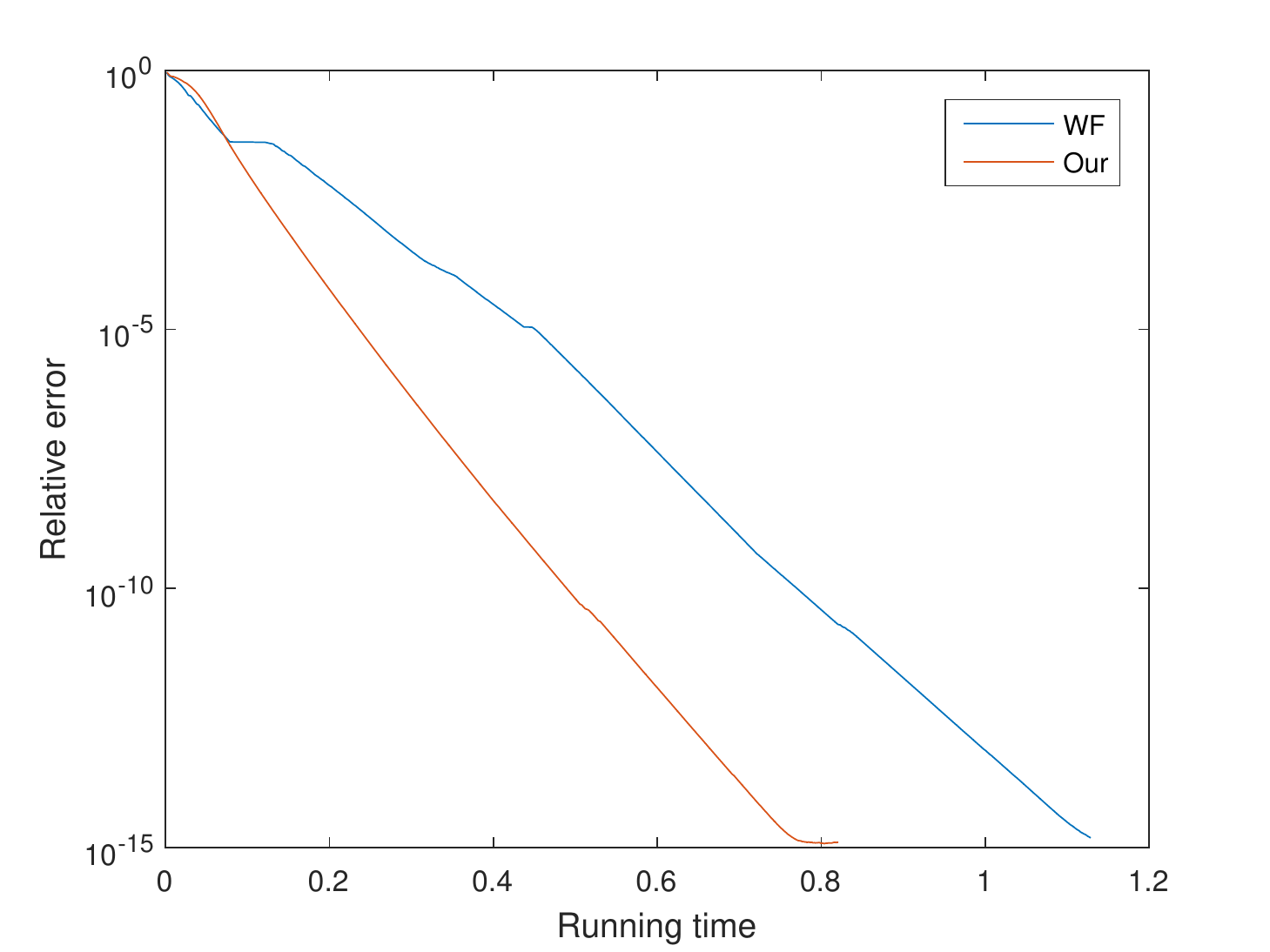}
\caption{The plot of relative error versus iteration number.}
\label{fig:iteration}
\end{center}
\end{figure}

Next, we demonstrate the efficiency of our proposed algorithm. We run WF 2500 times and our algorithm 1250 times for $\vx$ and $1250$ times for $\vy$, respectively.
 Figure \ref{fig:iteration} shows the plot of the relative error versus the iteration counts of our proposed method and the WF algorithm with $N=4.5d$ for the Gaussian model. From Figure \ref{fig:iteration}, we see that our propose algorithm can give $10^{-15}$ relative error for $\vx$ and $\vy$ after about $0.8$ seconds, but WF can not. This implies our algorithm needs less iterations than Wirtinger flow algorithm to get the same relative errors. 

\subsection{Real image data}
\label{real}
\begin{figure}[h]
\begin{center}
\includegraphics[width=0.90\textwidth]{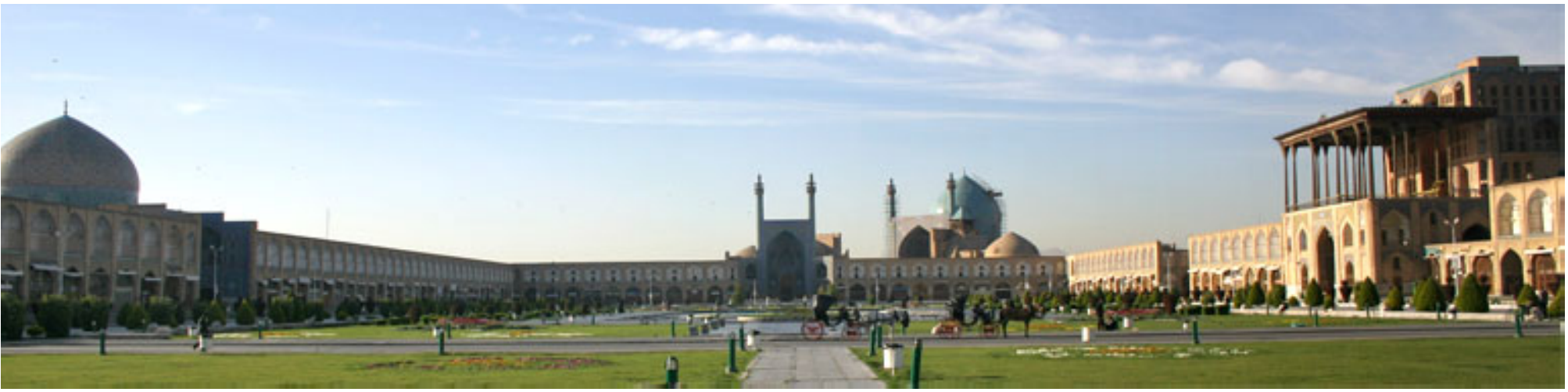}
\caption{The recovered images for Naqsh-e Jahan Square, Esfahan.}
\label{fig:naqsh}
\end{center}
\end{figure}

\begin{figure}[h]
\begin{center}
\includegraphics[width=0.9\textwidth]{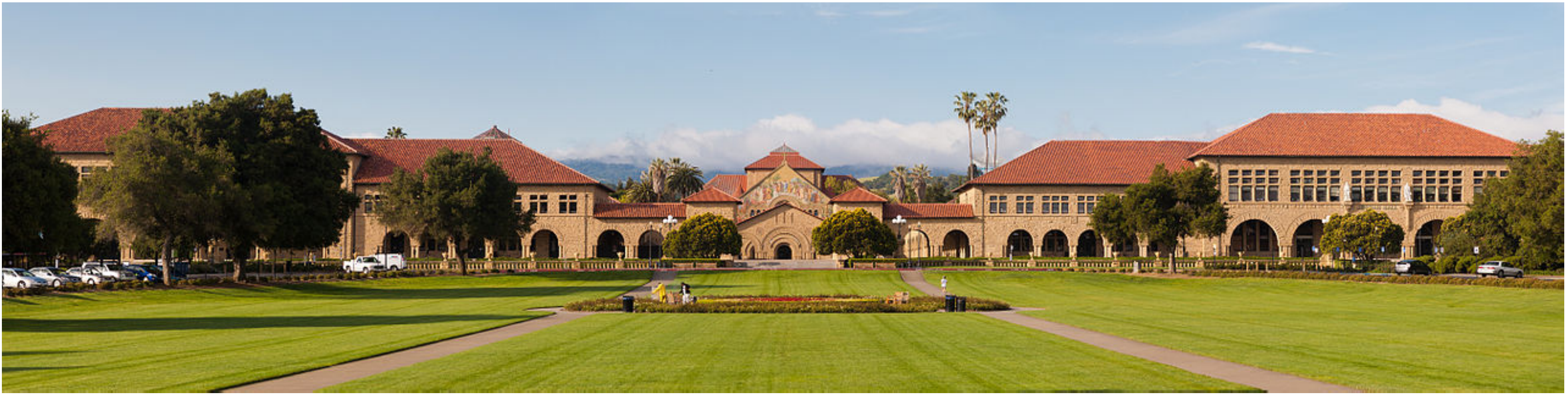}
\caption{The recovered images  for the stanford image.}
\label{fig:stanford}
\end{center}
\end{figure}

\begin{figure}[h]
\begin{center}
\includegraphics[width=0.45\textwidth]{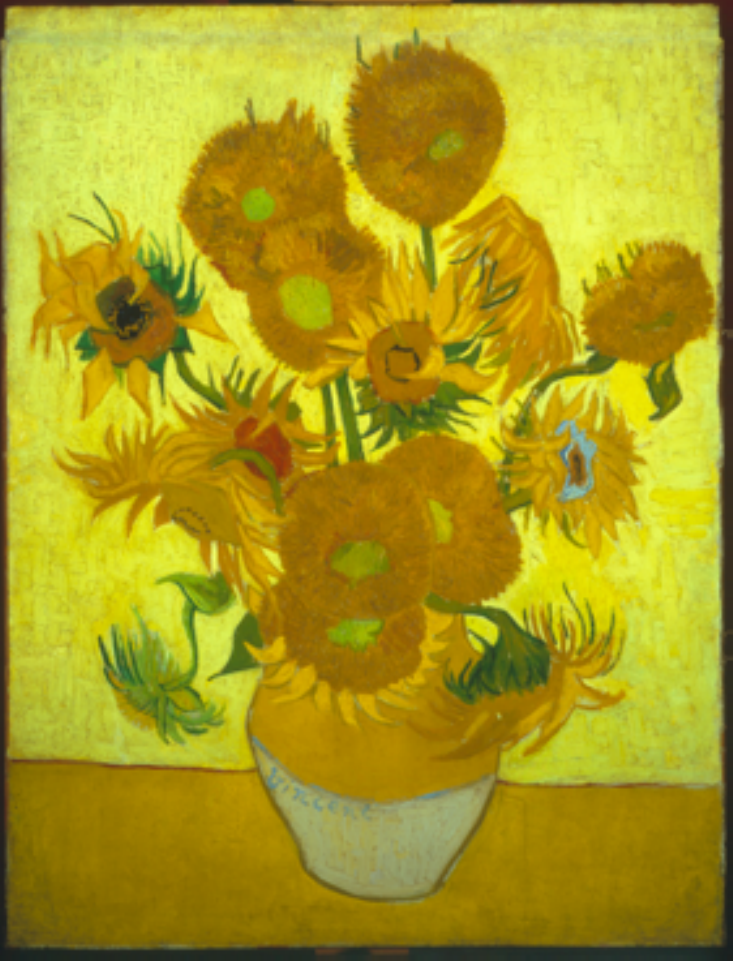}
\caption{The recovered images for the van Gogh painting $f458$.}
\label{fig:f458}
\end{center}
\end{figure}

We also test our algorithm for the CDF model on three real images in different sizes, namely, the Naqsh-e Jahan Square in the central Iranian city of Esfahan ($189\times 768\times 3$), the Stanford main quad ($320\times 1280\times 3$), and the van Gogh's painting $f458$ ($281\times 369\times 3$). Since those are all color images, we run our proposed algorithm and WF algorithm on each of the RGB channels. Let $\vx$ denote the underlying true image and $\hat{\vx}$ the solution by the algorithms. The relative error is defined as $\|\hat{\vx}-\vx\|/\|\vx\|$ with $\|\vx\|^2=\sum_{i,j,k}|x_{ijk}|^2$.  Table \ref{relative_error} lists the relative errors for WF with $2n$ iterations and our method with $n$ iterations for $\vx$ and $n$ iteration for $\vy$ with $n=100,125,150$. From the results in the table, we see that our proposed algorithm use less iterations than WF method to achieve the same relative error. In Figures \ref{fig:naqsh}, \ref{fig:stanford} and \ref{fig:f458}, the recoveries for the three real images are illustrated  after 150 iterations for $\vx$ and $150$ iterations for $\vy$.


\begin{table}
\begin{center}
\begin{tabular}{ c|ccc}
\hline\hline
\multicolumn{4}{ l }{ The Naqsh-e Jahan Square.} ($L=15$, $\tilde{\tau}_0=330$, $\tilde{\tau}_0=150$,\\
\multicolumn{4}{ l }{ $\tilde{\mu}_{max}=0.4$, $\mu_{max}=1$, $\lambda_0=8000$ and $\xi=0.001$.)}\\
\hline
Method &100&125&150\\ \hline
WF &$3.7097\times 10^{-4}$&6.1209$\times 10^{-7}$& 1.2521$\times 10^{-9}$\\
Our& $\mathbf{1.4927\times 10^{-4}}$&$\mathbf{5.9163\times 10^{-8}}$&$\mathbf{1.2817\times 10^{-11}}$\\\hline\hline
\multicolumn{4}{ l }{The Stanford main quad.} ($L=15$, $\tau_0=330$, $\tau_1=150$,\\
\multicolumn{4}{ l }{$\tilde{\mu}_{max}=0.4$, $\mu_{max}=1$, $\lambda_0=8000$ and $\xi=0.001$.)}\\
\hline
Method &100&125&150\\ \hline
WF &0.5608&$9.9557\times 10^{-4}$& $1.4987\times 10^{-6}$\\
Our& $\mathbf{0.5310}$&$\mathbf{1.4925\times 10^{-4}}$&$\mathbf{3.6023\times 10^{-8}}$\\
\hline\hline
\multicolumn{4}{ l }{The van Gogh's painting $f458$.} ($L=15$, $\tilde{\tau}_0=330$, $\tau_0=100$,\\
\multicolumn{4}{ l }{ $\tilde{\mu}_{max}=0.4$, $\mu_{max}=0.5$, $\lambda_0=5000$, $\xi=0.0015$.)}\\
\hline
Method & 100&125&150\\ \hline
WF &$0.2178$& $0.0028$ & $3.2730\times 10^{-6}$\\
Our& $\mathbf{7.7887\times 10^{-4}}$&$\mathbf{1.6263\times 10^{-6}}$&$\mathbf{2.6466\times 10^{-8}}$\\
\hline\hline
\end{tabular}
\caption{The relative errors.}
\label{relative_error}
\end{center}
\end{table}


\section*{Acknowledgement}
The authors would like to thank Emmanuel Candes, Mo Mu and Aditya Viswanathan for very helpful discussions.


{\small
\bibliographystyle{siam}
\bibliography{phase_retrieval}

\begin{thebibliography}{10}

\bibitem{alexeev2014phase}
{\sc Boris Alexeev, Afonso~S. Bandeira, Matthew Fickus, and Dustin~G. Mixon},
  {\em Phase retrieval with polarization}, SIAM Journal on Imaging Sciences, 7
  (2014), pp.~35--66.

\bibitem{balan2009painless}
{\sc Radu Balan, Bernhard~G. Bodmann, Peter~G. Casazza, and Dan Edidin}, {\em
  Painless reconstruction from magnitudes of frame coefficients}, Journal of
  Fourier Analysis and Applications, 15 (2009), pp.~488--501.

\bibitem{balan2006signal}
{\sc Radu Balan, Pete Casazza, and Dan Edidin}, {\em On signal reconstruction
  without phase}, Applied and Computational Harmonic Analysis, 20 (2006),
  pp.~345--356.

\bibitem{bolte2014proximal}
{\sc J{\'e}r{\^o}me Bolte, Shoham Sabach, and Marc Teboulle}, {\em Proximal
  alternating linearized minimization for nonconvex and nonsmooth problems},
  Mathematical Programming, 146 (2014), pp.~459--494.

\bibitem{candes2015phase1}
{\sc Emmanuel~J. Cand{\`e}s, Yonina~C. Eldar, Thomas Strohmer, and Vladislav
  Voroninski}, {\em Phase retrieval via matrix completion}, SIAM review, 57
  (2015), pp.~225--251.

\bibitem{candes2014solving}
{\sc Emmanuel~J. Cand{\`e}s and Xiaodong Li}, {\em Solving quadratic equations
  via phaselift when there are about as many equations as unknowns},
  Foundations of Computational Mathematics, 14 (2014), pp.~1017--1026.

\bibitem{candes2015phase}
{\sc Emmanuel~J. Cand{\`e}s, Xiaodong Li, and Mahdi Soltanolkotabi}, {\em Phase
  retrieval via {W}irtinger flow: Theory and algorithms}, Information Theory,
  IEEE Transactions on, 61 (2015), pp.~1985--2007.

\bibitem{candes2013phaselift}
{\sc Emmanuel~J. Cand{\`e}s, Thomas Strohmer, and Vladislav Voroninski}, {\em
  Phaselift: Exact and stable signal recovery from magnitude measurements via
  convex programming}, Communications on Pure and Applied Mathematics, 66
  (2013), pp.~1241--1274.

\bibitem{chai2010array}
{\sc Anwei Chai, Miguel Moscoso, and George Papanicolaou}, {\em Array imaging
  using intensity-only measurements}, Inverse Problems, 27 (2010), p.~015005.

\bibitem{chen2015solving}
{\sc Yuxin Chen and Emmanuel Candes}, {\em Solving random quadratic systems of
  equations is nearly as easy as solving linear systems}, in Advances in Neural
  Information Processing Systems, 2015, pp.~739--747.

\bibitem{hardt2014understanding}
{\sc Moritz Hardt}, {\em Understanding alternating minimization for matrix
  completion}, in Foundations of Computer Science (FOCS), 2014 IEEE 55th Annual
  Symposium on, IEEE, 2014, pp.~651--660.

\bibitem{harrison1993phase}
{\sc Robert~W. Harrison}, {\em Phase problem in crystallography}, JOSA A, 10
  (1993), pp.~1046--1055.

\bibitem{heinosaari2013quantum}
{\sc Teiko Heinosaari, Luca Mazzarella, and Michael~M. Wolf}, {\em Quantum
  tomography under prior information}, Communications in Mathematical Physics,
  318 (2013), pp.~355--374.

\bibitem{iwen2015robust}
{\sc Mark Iwen, Aditya Viswanathan, and Yang Wang}, {\em Robust sparse phase
  retrieval made easy}, Applied and Computational Harmonic Analysis,  (2015).

\bibitem{jain2013low}
{\sc Prateek Jain, Praneeth Netrapalli, and Sujay Sanghavi}, {\em Low-rank
  matrix completion using alternating minimization}, in Proceedings of the
  forty-fifth annual ACM symposium on Theory of computing, ACM, 2013,
  pp.~665--674.

\bibitem{millane1990phase}
{\sc Rick~P. Millane}, {\em Phase retrieval in crystallography and optics},
  JOSA A, 7 (1990), pp.~394--411.

\bibitem{sun2016geometric}
{\sc Ju~Sun, Qing Qu, and John Wright}, {\em A geometric analysis of phase
  retrieval}, arXiv preprint arXiv:1602.06664,  (2016).

\bibitem{tanner2016low}
{\sc Jared Tanner and Ke~Wei}, {\em Low rank matrix completion by alternating
  steepest descent methods}, Applied and Computational Harmonic Analysis, 40
  (2016), pp.~417--429.

\bibitem{waldspurger2015phase}
{\sc Ir{\`e}ne Waldspurger, Alexandre d’Aspremont, and St{\'e}phane Mallat},
  {\em Phase recovery, maxcut and complex semidefinite programming},
  Mathematical Programming, 149 (2015), pp.~47--81.

\bibitem{walther1963question}
{\sc Adriaan Walther}, {\em The question of phase retrieval in optics}, Journal
  of Modern Optics, 10 (1963), pp.~41--49.

\bibitem{wei2016guarantees}
{\sc Ke~Wei, Jian-Feng Cai, Tony~F Chan, and Shingyu Leung}, {\em Guarantees of
  riemannian optimization for low rank matrix completion}, arXiv preprint
  arXiv:1603.06610,  (2016).

\bibitem{zhangs2017median}
{\sc Huishuai Zhang, Yuejie Chi, and Yingbin Liang}, {\em Median-truncated
  nonconvex approach for phase retrieval with outliers},  (2017).

\end{thebibliography}
}
\appendix
\section{Larger step size}
In this appendix, we demonstrate that, when $\vx$ and $\vy$ are sufficiently close, our alternating gradient descent algorithm is roughly 1.5 times faster than the WF algorithm in the real case. 

To this end, we let $E(\vx,\vy)$ be the function defined in (\ref{E_function}). Choose $\lambda=0$ and assume $\vx_k\approx \vy_k$. Then, by Taylor's expansion, we obtain
\begin{equation}\label{eq:dec}
\begin{split}
&E(\vx_{k+1},\vy_k)\\
\approx&
E(\vx_{k},\vy_k)+\nabla_{\vx}E(\vx_{k},\vy_k)^*\left[\begin{matrix}\vx_{k+1}-\vx_k\cr \overline{\vx_{k+1}-\vx_k}\end{matrix}\right]+\frac12
\left[\begin{matrix}\vx_{k+1}-\vx_k\cr \overline{\vx_{k+1}-\vx_k}\end{matrix}\right]^*\nabla_{\vx}^2E(\vx_k,\vy_k)
\left[\begin{matrix}\vx_{k+1}-\vx_k\cr \overline{\vx_{k+1}-\vx_k}\end{matrix}\right]\cr
=&E(\vx_{k},\vx_k)-\alpha_k\left\|\nabla_{\vx}E(\vx_{k},\vy_k)\right\|_2^2+\frac{\alpha^2_k}{2}\Re\left(\nabla_{\vx}E(\vx_k,\vy_k)^*\left(\sum_n\vf_n\vf^*_n\vy_k\vy_k^*\vf_n\vf^*_n\right)\nabla_{\vx}E(\vx_{k},\vy_k)\right)\cr
=&E(\vx_{k},\vy_k)-\alpha_k\left\|(\nabla_{\vx}E(\vx_{k},\vy_k)\right\|_2^2
+\frac{\alpha_k^2}{2}\nabla_{\vx}E(\vx_{k},\vy_k)^*\left(\sum_n\vf_n\vf^*_n\vy_k\vy_k^*\vf_n\vf^*_n\right)\nabla_{\vx}E(\vx_{k},\vy_k)\cr
=&E(\vx_{k},\vy_k)-\alpha_k\left(\left\|\nabla_{\vx}E(\vx_{k},\vy_k)\right\|_2^2
-\frac{\alpha_k}{2}\nabla_{\vx}E(\vx_{k},\vy_k)^*\left(\sum_n\vf_n\vf^*_n\vy_k\vy_k^*\vf_n\vf^*_n\right)\nabla_{\vx}E(\vx_{k},\vy_k)\right)
\end{split}
\end{equation}
The last equality hold because $\sum_n\vf_n\vf^*_n\vy_k\vy_k^*\vf_n\vf^*_n$ is Hermitian. 
We choose $\alpha_k>0$. Therefore, $E(\vx_{k+1},\vy_k)-E(\vx_{k},\vy_k)\leq 0$ as long as
$$
\frac{2}{\alpha_k}\geq\frac{\nabla_{\vx}E(\vx_{k},\vy_k)^*\left(\sum_n\vf_n\vf^*_n\vy_k\vy_k^*\vf_n\vf^*_n\right)\nabla_{\vx}E(\vx_{k},\vy_k)}{\left\|\nabla_{\vx}E(\vx_{k},\vy_k)\right\|_2^2},
$$
which is guaranteed if 
$$
\alpha_k\leq \frac{2}{\|\sum_n\vf_n\vf^*_n\vy_k\vy_k^*\vf_n\vf^*_n\|_2}.
$$
To minimize $E(\vx_{k+1},\vy_k)-E(\vx_{k},\vy_k)$, it is easy seen from \eqref{eq:dec} that $\alpha_k$ is chosen as
$$
\alpha_k=\frac{\left\|\nabla_{\vx}E(\vx_{k},\vy_k)\right\|_2^2}{\nabla_{\vx}E(\vx_{k},\vy_k)^*\left(\sum_n\vf_n\vf^*_n\vy_k\vy_k^*\vf_n\vf^*_n\right)\nabla_{\vx}E(\vx_{k},\vy_k)}.
$$
In this case,
\begin{equation}\label{eq:decfx}
E(\vx_{k+1},\vy_k)-E(\vx_{k},\vy_k)
\approx
-\frac{\left\|\nabla_{\vx}E(\vx_{k},\vy_k)\right\|_2^4}{2\nabla_{\vx}E(\vx_{k},\vy_k)^*\left(\sum_n\vf_n\vf^*_n\vy_k\vy_k^*\vf_n\vf^*_n\right)\nabla_{\vx}E(\vx_{k},\vy_k)}.
\end{equation}

Now we consider the WF algorithm, which minimizes $G(\vx)=\frac{1}{N}\sum^N_{n=1}(|\vf^*_n\vx|^2-b_n)^2$. Assume we have the same $\vx_k$ as in the alternating gradient descent algorithm, and the WF algorithm generates the new iterates by
$$
\vx_{k+1}=\vx_k-\delta_k\nabla_{\vx}G(\vx_k).
$$
With the optimal choice of $\delta_k$, an analogous analysis leads to 
\begin{equation}\label{eq:decg}
G(\vx_{k+1})-G(\vx_{k})
\approx
-\frac12\frac{\left\|\nabla_{\vx}G(\vx_k)\right\|_2^4}{\Re\left(\nabla_{\vx}G(\vx_k)^*H_{11}(\vx_k)\nabla_{\vx}G(\vx_k)\right)+\Re\left(\nabla_{\vx}G(\vx_k)^TH_{21}(\vx_k)\nabla_{\vx}G(\vx_k)\right)},
\end{equation}
where $\Re$ denotes the real part, and
\begin{eqnarray*}
&H_{11}(\vx_k)
=4\sum_n\vf_n\vf^*_n\vx_k\vx_k^*\vf_n\vf^*_n,\\
&H_{21}(\vx_k)=\sum_n\left(\overline{\vf_n\vf^*_n\vx_k}\vx_k^*\vf_n\vf^*_n+\overline{\vf_n\vf^*_n\vx_k}\vx_k^*\vf_n\vf^*_n\right).
\end{eqnarray*}

%


Since we assumed $\vx_k\approx\vy_k$ and $\lambda=0$,
\begin{equation}\label{eq:gra}
\nabla_{\vx}G(\vx_k)\approx2\nabla_{\vy}E(\vx_k,\vy_k),
\end{equation}
which implies
\begin{equation}\label{eq:den1}
\begin{split}
\Re\left(\nabla_{\vx}G(\vx_k)^*H_{11}(\vx_k)\nabla_{\vx}G(\vx_k)\right)
=&\nabla_{\vx}G(\vx_k)^*H_{11}(\vx_k)\nabla_{\vx}G(\vx_k)\cr
\approx&\left(2\nabla_{\vy}E(\vx_k,\vy_k)\right)^*\left(4\sum_n\vf_n\vf^*_n\vx_k\vx_k^*\vf_n\vf^*_n\right)\left(2\nabla_{\vy}E(\vx_k,\vy_k)\right)\cr
\approx&16\cdot\left(\nabla_{\vy}E(\vx_k,\vy_k)\right)^*\left(\sum_n\vf_n\vf^*_n\vy_k\vy_k^*\vf_n\vf^*_n\right)\nabla_{\vy}E(\vx_k,\vy_k).
\end{split}
\end{equation}
If we further assume all vectors involved are real, then we have $H_{21}(\vx_k)=2\sum_n\vf_n\vf^*_n\vx_k\vx_k^*\vf_n\vf^*_n$ and
\begin{equation}\label{eq:den2}
\begin{split}
\Re\left(\nabla_{\vx}G(\vx_k)^TH_{21}(\vx_k)\nabla_{\vx}G(\vx_k)\right)
&=\nabla_{\vx}G(\vx_k)^*H_{21}(\vx_k)\nabla_{\vx}G(\vx_k)\cr
&\approx8\cdot\nabla_{\vy}E(\vx_k,\vy_k)^*\left(\sum_n\vf_n\vf^*_n\vy_k\vy_k^*\vf_n\vf^*_n\right)\nabla_{\vy}E(\vx_k,\vy_k).
\end{split}
\end{equation}
Substituting \eqref{eq:gra}, \eqref{eq:den1}, and \eqref{eq:den2} into \eqref{eq:decg}, we get
$$
\left(G(\vx_{k+1})-G(\vx_k)\right)
\approx
\frac23\left(E(\vx_{k+1},\vy_k)-E(\vx_k,\vy_k)\right)
$$
This means that \emph{the alternating gradient descent algorithm is $1.5$ times faster than Wirtinger flow in terms of the decreasing of the objective.}

\end{document}